\documentclass[12pt]{article}

\usepackage{amsmath,amssymb}
\begin{document}

\pagestyle{myheadings} \markright{MELLIN TRANSFORMS OF WHITTAKER
FUNCTIONS}

\title{Mellin transforms of Whittaker functions}
\author{Anton Deitmar}

\date{}
\maketitle

$$ $$

{\bf Abstract.} In this note we show that for an arbitrary
semisimple Lie group and any admissible irreducible Banach
representation the Mellin transforms of Whittaker functions
extend to meromorphic functions. We determine the possible poles.

\tableofcontents

\def \1{{\bf 1}}
\def \a{{{\mathfrak a}}}
\def \ad{{\rm ad}}
\def \al{\alpha}
\def \ar{{\alpha_r}}
\def \A{{\mathbb A}}
\def \Ad{{\rm Ad}}
\def \Aut{{\rm Aut}}
\def \b{{{\mathfrak b}}}
\def \bs{\backslash}
\def \B{{\cal B}}
\def \c{{\mathfrak c}}
\def \cent{{\rm cent}}
\def \C{{\mathbb C}}
\def \CA{{\cal A}}
\def \CB{{\cal B}}
\def \CC{{\cal C}}
\def \CD{{\cal D}}
\def \CE{{\cal E}}
\def \CF{{\cal F}}
\def \CG{{\cal G}}
\def \CH{{\cal H}}
\def \CHC{{\cal HC}}
\def \CL{{\cal L}}
\def \CM{{\cal M}}
\def \CN{{\cal N}}
\def \CP{{\cal P}}
\def \CQ{{\cal Q}}
\def \CO{{\cal O}}
\def \CS{{\cal S}}
\def \CT{{\cal T}}
\def \CV{{\cal V}}
\def \CW{{\cal W}}
\def \det{{\rm det}}
\def \diag{{\rm diag}}
\def \dist{{\rm dist}}
\def \End{{\rm End}}
\def \eqn{\begin{eqnarray*}}
\def \endeqn{\end{eqnarray*}}
\def \F{{\mathbb F}}
\def \Fx{{\mathfrak x}}
\def \FX{{\mathfrak X}}
\def \g{{{\mathfrak g}}}
\def \ga{\gamma}
\def \Ga{\Gamma}
\def \Gal{{\rm Gal}}
\def \h{{{\mathfrak h}}}
\def \Hom{{\rm Hom}}
\def \im{{\rm im}}
\def \Im{{\rm Im}}
\def \Ind{{\rm Ind}}
\def \k{{{\mathfrak k}}}
\def \K{{\cal K}}
\def \l{{\mathfrak l}}
\def \la{\lambda}
\def \lap{\triangle}
\def \li{{\rm li}}
\def \La{\Lambda}
\def \Lie{{\rm Lie}}
\def \m{{{\mathfrak m}}}
\def \mod{{\rm mod}}
\def \n{{{\mathfrak n}}}
\def \name{\bf}
\def \Mat{{\rm Mat}}
\def \N{\mathbb N}
\def \o{{\mathfrak o}}
\def \ord{{\rm ord}}
\def \O{{\cal O}}
\def \p{{{\mathfrak p}}}
\def \ph{\varphi}
\def \prf{\noindent{\bf Proof: }}
\def \Per{{\rm Per}}
\def \q{{\mathfrak q}}
\def \qed{$ $\newline $\frac{}{}$\hfill {\rm Q.E.D.}\vspace{15pt}\pagebreak[0]}
\def \Q{\mathbb Q}
\def \res{{\rm res}}
\def \R{{\mathbb R}}
\def \Re{{\rm Re \hspace{1pt}}}
\def \r{{\mathfrak r}}
\def \ra{\rightarrow}
\def \rank{{\rm rank}}
\def \rop{\rho_{_P}}
\def \supp{{\rm supp}}
\def \Spin{{\rm Spin}}
\def \t{{{\mathfrak t}}}
\def \T{{\mathbb T}}
\def \tr{{\hspace{1pt}\rm tr\hspace{2pt}}}
\def \vol{{\rm vol}}
\def \z{\zeta}
\def \Z{\mathbb Z}
\def \={\ =\ }

\newcommand{\rez}[1]{\frac{1}{#1}}
\newcommand{\der}[1]{\frac{\partial}{\partial #1}}
\newcommand{\norm}[1]{\parallel #1 \parallel}
\renewcommand{\matrix}[4]{\left( \begin{array}{cc}#1 & #2 \\ #3 & #4 \end{array}
            \right)}
\renewcommand{\sp}[2]{\langle #1,#2\rangle}

\newcounter{lemma}
\newcounter{corollary}
\newcounter{proposition}
\newcounter{theorex}

\newtheorem{conjecture}{\hspace{-20pt}\stepcounter{lemma} \stepcounter{corollary}
    \stepcounter{proposition}\stepcounter{theorex}Conjecture}[section]
\newtheorem{lemma}{\hspace{-20pt}\stepcounter{conjecture}\stepcounter{corollary}
    \stepcounter{proposition}\stepcounter{theorex}Lemma}[section]
\newtheorem{corollary}{\hspace{-20pt}\stepcounter{conjecture}\stepcounter{lemma}
    \stepcounter{proposition}\stepcounter{theorex}Corollary}[section]
\newtheorem{proposition}{\hspace{-20pt}\stepcounter{conjecture}\stepcounter{lemma}
    \stepcounter{corollary}\stepcounter{theorex}Proposition}[section]

\newtheorem{theorex}{\hspace{-40pt}\stepcounter{conjecture} \stepcounter{lemma}
    \stepcounter{corollary} \stepcounter{proposition}Theorem}[section]
\newenvironment{theorem}{\vspace{30pt}\begin{theorex}}{\end{theorex}\vspace{30pt}}

\newpage

\section{Whittaker functions}
Let $G$ be a semisimple connected Lie group with finite center.
Fix a maximal compact subgroup $K$ and let $\theta$ denote the
Cartan involution with fixed point set $K$. We will write $\g_0$
for the real Lie algebra of $G$ and $\g$ for its
complexification. Next $U(\g)$ will denote the universal
enveloping algebra of $\g$. We will interpret $U(\g)$ as the
algebra of all left invariant differential operators on $G$.

Let $P=MAN$ be a minimal parabolic subgroup, where we assume that
$A$ and hence also $M$ are $\theta$-stable.

Let $Z(A)$ be the centralizer of $A$ in $G$ and $N(A)$ be its
normalizer. The Weyl group
$$
W\= W(A,G)\= N(A)/Z(A)
$$
is a finite group. The Bruhat decomposition
$$
G\= \bigcup_{w\in W} PwP\= \bigcup_{w\in W} PwN
$$
is a disjoint decomposition.

Let $\m_0,\a_0,\n_0$ be the real Lie algebras of the groups $M, A,
N$. and let $\m,\a,\n$ be their complexifications. Let
$\Phi^+=\Phi^+(\a,\g)$ be the set of positive roots of the pair
$(\a,\g)$ given by the choice of the parabolic $P$. For
$\alpha\in \Phi^+(\a,\g)$ let $\g_\alpha$ denote its root space
and let $m_\alpha=\dim\g_\alpha$. Further let $\n_\alpha
=\g_\alpha\cap\g_0$. Then
$\n_0=\bigoplus_{\alpha\in\Phi^+}\n_\alpha$. Let $\rho_{_P}=\rez
2\sum_{\alpha\in\Phi^+} m_\alpha \alpha$ be the modular shift of
$P$.

The Killing form on $\g$ is positive definite on $\a_0$. It
induces an identification of $\a_0$ with its dual $\a_0^*$ and
also a bilinear form $\langle .,.\rangle$ on $\a_0^*$.

Let $\pi$ be a continuous representation of $G$ on some Banach
space. Let $\pi_\infty$ the Fr\'echet space of differentiable
vectors for $\pi$. Fix a continuous linear functional
$\psi=\psi_\eta$ on $\pi_\infty$ such that
$$
\psi(\pi(n)\phi)\= \eta(n)\psi(\phi)
$$
for every $\phi\in \pi_\infty$. Such a $\psi$ is called a {\it
Whittaker functional} to the character $\eta$. For $\phi\in
\pi_\infty$ set
$$
W_\phi(x)\=\psi(\pi(x)\phi),\ \ \ x\in G,
$$
the corresponding {\it Whittaker function} on the group $G$.

\begin{proposition}\label{2.1} Assume that $\pi$ is admissible and
of finite length. There is a finite set $\{ D_j\}\subset U(\g)$
and a natural number $k_0$, depending on $\pi$, such that for all
$\phi\in \pi_\infty$ we have
$$
|W_\phi(am)|\ \le\ \sum_{w\in W}a^{k_0w\rho_{_P}}\ \sum_j
\norm{D_j\phi}.
$$
\end{proposition}

\prf The continuity of $\psi$ implies that there is a finite set
$\{ D_j'\}\subset U(\g)$ such that
$$
|W_\phi(1)|\le \sum_j \norm{D_j'\phi}.
$$
From \cite{Wall-asymp}, Lemma 2.2 we derive the existence of
$c>0$, $k'\in\N$ such that the operator norm can be estimated:
$$
\norm{\pi(am)}\ \le\ c\sum_{w\in W} a^{k'w\rho_{_P}}
$$
for $a\in A$, $m\in M$. We thus get
\eqn
W_\phi(am) &=& W_{\pi(am)\phi}(1)\\
&\le&\sum_j\norm{D_j'\pi(am)\phi}\\
&=&\sum_j\norm{\pi(am)\Ad(am)^{-1}D_j'\phi}\\
&\le& c\sum_{w\in W}a^{k'w\rho_{_P}}\sum_j \norm{\Ad(am)^{-1}D_j'\phi}.
\endeqn
Let $U(\g)^\nu$ be the finite dimensional space of all $D\in
U(\g)$ of degree $\le\nu$. Let $\{D_j''\}$ be a basis of
$U(\g)^\nu$. For $\nu$ large enough we get
$$
\Ad(am)^{-1} D_j'\= \sum_i a_{i,j}(am)D_i''.
$$
By the properties of the adjoint action it follows that there is
a constant $c_1>0$ with $a_{i,j}(am)\le c_1\sum_{w\in
W}a^{2\nu\rho_{_P}}$. The lemma follows with $D_j$ being a
multiple of $D_j''$. \qed

From the preceding lemma we will now conclude that the Whittaker
function $W_\phi(am)$ actually is rapidly decreasing.

\begin{proposition} \label{1.2} Let $\pi$ be admissible and of finite length.
Let $\alpha_0\in\Phi^+(\a,\g)$ be a positive root such that
$\log\eta$ is nontrivial on the root space
$\n_{\alpha_0}\subset\n$. For every natural number $N$ there are
$D_1,\dots,D_m\in U(\g)$ such that for every $\phi\in \pi_\infty$
and every $a\in A$ we have
$$
|W_\phi(am)|\ \le\ a^{-N\alpha_0}\left(\sum_{w\in
W}a^{k_0w\rho_{_P}}\right)\left(\sum_{j=1}^m
\norm{D_j\phi}\right).
$$
\end{proposition}

\prf Let $X_1,\dots,X_n$ be a basis of the root space
$\n_{\alpha_0}$. Since $\eta$ is nontrivial on $\n_{\alpha_0}$
the function
$$
f(m)\=\sum_j |\log\eta(\Ad(m)X_j)|
$$ is nowhere vanishing on $M$. Since $M$ is compact there is
$c>0$ with $f(m)\ge c$ for all $m\in M$. It follows
\eqn
W_{X_j\phi}(am) &=& \frac{d}{dt} W_\phi(am\exp(tX_j))|_{t=0}\\
&=& \frac{d}{dt} W_\phi(\exp(\Ad(am)tX_j)am)|_{t=0}\\
&=& \frac{d}{dt}\eta(\exp(t\Ad(am)X_j))|_{t=0} W_\phi(am)\\
&=& \log\eta(\Ad(am)X_j)W_\phi(am)\\
&=& a^{\alpha_0}\log\eta(\Ad(m)X_j)W_\phi(am).
\endeqn
So for $a\in A$ it follows
\eqn
\sum_j|W_{X_j\phi}(am)| &=&
|W_\phi(am)|a^{\alpha_0}\sum_j|\log\eta(\Ad(m)X_j)|\\
&\ge& c a^{\alpha_0} |W_\phi(am)|.
\endeqn
Iterating this process gives for arbitrary $N\in\N$ the existence
of $\{ D_j'\}\subset U(\n_{\alpha_0})$ such that
$$
|W_\phi(am)|\ \le\ a^{-N\alpha_0}\sum_j|W_{D_j'\phi}(am)|.
$$
Applying the last lemma to $D_j'\phi$ gives the claim.
\qed

A character $N\ra\T$ factors over $N^{ab}=N/[N,N]$, where $[N,N]$
denotes the closed subgroup of $N$ generated by all commutators
$aba^{-1}b^{-1}$ for $a,b\in N$. Then $N^{ab}$ is an abelian,
simply connected Lie group with Lie algebra $\n^{ab}=\n/[\n,\n]$,
where in this case $[.,.]$ denotes the Lie bracket. It follows
that $N^{ab}$ is isomorphic to its Lie algebra and the characters
of $N$ thus identify with the linear functionals on $\n^{ab}$. Let
$$
\n_{simp}\= \bigoplus_{^{\alpha\in\Phi^+(\a,\g)}_{\alpha\ \rm
simple}} \n_{\alpha}.
$$

\begin{lemma}
We have
$$
\n_0\=\n_{0,simp}\oplus [\n_0,\n_0].
$$
\end{lemma}

\prf This follows from Proposition 8.4 d) of \cite{humpliealg}.
\qed

Thus each linear functional on $\n_{simp}$ extends to a character
on $N$. This implies that there are characters $\eta$ which are
nontrivial on $\n_{\alpha}$ for each simple root $\alpha$. In
this case we call $\eta$ a {\it generic} character. Let
$\Delta\subset\Phi^+$ be the set of simple roots.

\begin{corollary} \label{1.4}
Suppose that $\pi$ is admissible of finite length and that the
character $\eta$ is generic. For every natural number $N$ there
are $D_1,\dots,D_m\in U(\g)$ such that for every $\phi\in
\pi_\infty$ and every $a\in A$ we have
$$
|W_\phi(am)|\ \le\ \min_{\alpha\in\Delta}
a^{-N\alpha}\left(\sum_{w\in
W}a^{k_0w\rho_{_P}}\right)\left(\sum_{j=1}^m
\norm{D_j\phi}\right).
$$
\end{corollary}

\prf This is a direct consequence of the proposition.
\qed

To be able to apply this corollary we have to make sure that
there are generic characters which are trivial on $\Ga\cap N$. By
\cite{rag}, Theorem 1.13 we infer that the image of $\Ga$ in
$N^{ab}$ is a lattice, which implies that $N^{ab}/\Ga$ is a
torus. This implies the existence of an abundance of generic
characters that are trivial on $\Ga$.

\section{The Mellin transform, rank one}\label{rank-one}
Let $\tau$ be a finite dimensional unitary representation of the
compact group $M$. Let $\xi_{\tau}$ be an arbitrary matrix
coefficient of the representation $\tau$. For $\phi\in
\pi_\infty$ and $\la\in \a^*$ let
$$
I_\phi(\xi_{\tau},\la)\= \int_A\int_M
W_\phi(am)\overline{\xi_{\tau}(m)} a^{\la-\rho_{_P}} da dm.
$$
For $\mu\in\a_0^*$ we write $\mu>0$ if $\sp{\mu}{\alpha}>0$ for
all $\alpha\in\Phi^+$ and $\mu>\nu$ if $\mu-\nu>0$. For the rest
of this section we assume that $\pi$ is admissible of finite
length. From the last two lemmas we get

\begin{proposition}\label{2.4}
Suppose $\eta$ is generic, then the integral
$I_\phi(\xi_{\tau},\la)$ converges absolutely for $\Re(\la)
>k_0\rho_{_P}$, where $k_0$ is the number of
Proposition \ref{2.1}. The linear functional $\phi\mapsto
I_\phi(\xi_{\tau},\la)$ is continuous on $\pi_\infty$.
\end{proposition}

Let $(X_j)$ be a basis of $\n$ such that each
$X_j\in\n_{\alpha_j}$ for some root $\alpha_j$. Write $Z(\g)$ for
the center of $U(\g)$.

\begin{lemma}
Let $D\in Z(\g)$, then $D$ can be written as
$$
D=D_{AM}+\sum_j X_jD_j
$$
with $D_{AM}\in Z(\a\oplus\m)$ and $D_j\in U(\g)$. Moreover we
have $D_{AM}\in Z(\a\oplus \m)^M$ the subalgebra of
$M$-invariants.
\end{lemma}

\prf Write $\n^-=\theta(n)$, then we have the decomposition
$$
\g\= \n^-\oplus \a\oplus\m\oplus\n
$$
and hence, by the Poincar\'e-Birkhoff-Witt Theorem:
$$
U(\g)\= U(\a\oplus\m)\oplus (\n U(\g)+U(\g)\n^-).
$$
So $D\in Z(\g)$ can be written as $D_{AM}+f$ with $f\in\n
U(\g)+U(\g)\n^-$. Since $G$ is connected we get $Z(\g)=U(\g)^G$,
so for $m\in M$ and $D\in Z(\g)$ we have $\Ad(m)D=D$. The
decomposition above is stable under $M$, so it follows that
$\Ad(m)D_{AM}=D_{AM}$, which implies
$$
Z(\g)\ \subset\ Z(\a\oplus \m)^M\oplus(\n U(\g)+U(\g)\n^-).
$$
Writing $X_j^-=\theta(X_j)$ we see by the
Poincar\'e-Birkhoff-Witt Theorem that $f$ is a sum of monomials
of the form
$$
X^aD_1(X^-)^b,
$$
where for $a,b\in \Z_+^d$ with $d=\dim N$:
$$
X^a\= X_1^{a_1}\dots X_d^{a_d},\ \ \ (X^-)^b\= (X_1^-)^{b_1}\dots
(X_d^-)^{b_d},
$$
and $D_1\in U(\a\oplus\m)$. For $H\in \a$ we compute that
$$
[H,X^aD_1(X^-)^b]
$$
equals
$$
(a_1\alpha_1(H)+\dots +a_d\alpha_d(H)-b_1\theta(\alpha_1)(H)-\dots
-b_d\theta(\alpha_d)(H)) X^aD_1(X^-)^b.
$$
Since $f$ lies in $\n U(\g)+U(\g)\n^-$ it follows that only
monomials occur with $a$ and $b$ not both zero. Next, since $f$
commutes with each $H\in\a$ it follows that for each monomial
both $a$ and $b$ are nonzero, which implies the lemma.
\qed

Assume now that the representation $\pi$ is {\it quasi-simple},
which means that the center $Z(\g)$ acts by scalars. Let
$\wedge_\pi : Z(\g)\ra\C$ denote the {\it infinitesimal
character} of $\pi$, i.e.
$$
\pi(D)\phi\=\wedge_\pi(D)\phi
$$
holds for every $\phi\in \pi_\infty$, $D\in Z(\g)$. For
$\sigma\in\hat{M}$ and $\nu\in\a^*$ the representation
$$
\pi_{\sigma ,\nu}\= {\rm Ind}_P^G(\sigma\otimes
(\nu+\rho_{_P})\otimes 1)
$$
is known to be quasi-simple. Let $\wedge_{\sigma,\nu}$ denote its
infinitesimal character.

\begin{proposition}\label{2.6}
For $D\in Z(\g)$ write $D=D_{AM}+\sum_{j=1}^d X_j D_j$ as in the
last lemma. Let $r_j$ denote the matrix coefficient
$r_j(m)=\overline{\eta(-\Ad(m)X_j)}$. Then, for
$\langle\Re(\la),\rho_{_P}\rangle >0$ we have
$$
I_\phi(\xi_{\tau},\la)\=\frac{\sum_{j=1}^d
I_{D_j\phi}(r_j\xi_{\tau},\la+\alpha_j)}
{\wedge_\pi(D)-\wedge_{\tau,\rho_{_P}-\la}(D)}.
$$
\end{proposition}

\prf For $D\in Z(\g)$ we have on the one hand
$$
I_{D\phi}(\xi_{\tau},\la)\=\wedge_\pi(D)I_\phi(\xi_{\tau},\la),
$$
and on the other
$$
I_{D\phi}(\xi_{\tau},\la)\=
I_{D_{AM}\phi}(\xi_{\tau},\la)+\sum_{j=1}^d
I_{X_jD_j\phi}(\xi_{\tau},\la).
$$
We compute
\eqn
I_{X_j\phi}(\xi_{\tau},\la) &=& \int_A\int_M
W_{X_j\phi}(am)\overline{\xi_{\tau}(m)} a^{\la-\rho_{_P}} dm da\\
&=& \int_A\int_M a^{\alpha_j} \log\eta(\Ad(m)X_j)W_\phi(am)
\overline{\xi_{\tau}(m)} a^{\la-\rho_{_P}}dm da\\
&=& I_\phi(r_j\xi_{\tau},\la+\alpha_j)
\endeqn
To prove the proposition it remains to show
$$
I_{D_{AM}\phi}(\xi_{\tau},\la)\=\wedge_{\tau,\rho-\la}(D)
I_\phi(\xi_{\tau},\la).
$$
Fix $\la\in \a^*$ with $\Re(\la)$ large. Let $\breve{\tau}$ denote
the contragredient representation to $\tau$ and consider the
representation $\ga$ of $AM$ given by
$$
\ga\= (\la-\rho_{_P})\otimes \breve \tau.
$$
Note that the function $am\mapsto
\overline{\xi_{\tau}(m)}a^{\la-\rho_{_P}}$ is a matrix
coefficient of $\ga$, we denote it by $\xi_\ga(am)$.

On $U(\g)$ we have a unique $\C$-linear involutary
anti-autormorphism given by $X'=-X$ for $X\in\g$. For $D\in
U(\g)$ we have
$$
\int_G Df(x)g(x)dx\=\int_Gf(x)D'g(x)dx.
$$
We compute
\eqn
I_{D_{AM}\phi}(\xi_{\tau},\la) &=&\int_A\int_M
W_{D_{AM}\phi}(am)\overline{\xi_{\tau}(m)} a^{\la-\rho_{_P}} dm da\\
&=&\int_A\int_M
D_{AM}W_{\phi}(am)\xi_\ga(am) dm da\\
&=&\int_A\int_M
W_{\phi}(am)D_{AM}'\xi_\ga(am) dm da\\
&=& \wedge_{\tau,\rho_{_P}-\la}(D)I_\phi(\xi_{\tau},\la)
\endeqn
The proposition follows.
\qed

Let $\b\subset \m$ be a Cartan subalgebra and let $\CW
=W(\a\oplus \b,\g)$ be the big Weyl group. Via the Harish-Chandra
homomorphism the infinitesimal character $\wedge_\pi$ can be
viewed as (a Weyl group orbit of) an element of the dual space of
$\a\oplus \b$.

Let $r=\dim A$. If $r=1$ there are at most two positive roots in
$\Phi^+(\a,\g)$. Let $\alpha_0$ be the short positive root in
this case. Let $R$ denote the adjoint representation of $M$ on
$\n$. Let $\la_{\tau}\in\b^*$ be the infinitesimal character of
$\tau$.

\begin{theorem}\label{2.7}
If $r=1$ then the function $I_\phi(\xi_{\tau},\la)$ has a
meromorphic continuation to $\a^*$. There is a possible pole at
$\la$ only if the following conditions are satisfied: There are
integers $0\le l\le k$ and an irreducible subrepresentation $\ga$
of $\tau\otimes R^{\otimes k}$ such that, $\la_\ga$ denoting its
infinitesimal character we have that
$$
\la_\ga+\rho_{_P}-(\la+(k+l)\alpha_0)
$$
lies in the Weyl group orbit of $\wedge_\pi$.
\end{theorem}

\prf The function $I_\phi(\xi_{\tau},\la)$ is holomorphic in the
region $\Re(\la)>k_0\rho_{_P}$. Let $\alpha_1$ be the short
positive root. Now Proposition \ref{2.6} shows that
$I_\phi(\xi_{\tau},\la)$ extends to a meromorphic function on
$\Re(\la)>k_0\rho_{_P}-\alpha_1\}$ with possible poles where
$$
\wedge_\pi(D)\=\wedge_{\tau,\rho_{_P}-\la}(D)
$$
for every $D\in Z(\g)$. This can only be when the two
infinitesimal characters $\wedge_\pi$ and
$\wedge_{\tau,\rho_{_P}-\la}$ are in the same Weyl group orbit. We
iterate this replacing $\xi_{\tau}$ by $r_j\xi_{\tau}$ which is a
matrix coefficient of $\tau\otimes R$. Further iteration gives
the claim in the case $r=1$.
\qed

\section{The higher rank case}
For every simple root $\alpha\in\Delta$ let $G_\alpha$ be the
connected Lie subgroup of $G$ with Lie algebra generated by
$\n_{-\alpha}\oplus\n_\alpha$. Then $G_\alpha$ is semisimple of
real rank one with split torus
$$
A_\alpha\=\{ a\in A | a^\beta =1\ {\rm for\ all}\ \beta\in\Delta,
\beta\ne\alpha\}.
$$
Let $\a_{\alpha,0}$ be the Lie algebra of $A_\alpha$ and
$\a_\alpha$ its complexification. Let $M_\alpha =M\cap G_\alpha$
and $\tau_\alpha =\tau |_{M_\alpha}$.

The set $\Delta$ is a basis of $\a^*$, hence every $\la\in\a^*$
can uniquely be written as $\la = \sum_{\alpha\in\Delta}
\la_\alpha$, where $\la_\alpha\in\a_\alpha^*$. Let $m_\alpha$ be
the multiplicity of the root $\alpha$ and let $\rop'=\rop-\rez
2\sum_{\alpha\in\N\Delta}m_\alpha\alpha$.

\begin{theorem}\label{higher-rank} Let $\pi$ be an irreducible
admissible Banach representation of $G$. Let $\phi\in\pi^\infty$,
then the function $\la\mapsto I_\phi(\xi_\tau,\la)$ extends to a
meromorphic function on $\a^*$ with possible poles along the sets
$(\la+2\rop')_\alpha = c_\alpha$, where $\alpha\in\Delta$ and
$c_\alpha\in\a_\alpha^*$ is such that there is an irreducible
$G_\alpha$-subquotient $\pi_\alpha$ of $\pi$ and there are
integers $0\le l\le k$ and an irreducible subrepresentation $\ga$
of $\tau_\alpha\otimes (\Ad|\n_\alpha)^{\otimes k}$ such that,
$\wedge_\ga$ denoting its infinitesimal character we have that
$$
\wedge_\gamma+\rop^\alpha -(c_\alpha + (k+l)\alpha)
$$
lies in the Weyl orbit of $\wedge_{\pi_\alpha}$.
\end{theorem}

The Proof of this theorem will occupy the rest of the section. We
start by considering a special case. So let $(\sigma,V_\sigma)$
be an irreducible unitary representation of $M$. Since $M$ is
compact it follows that $\sigma$ is finite dimensional. Let
$\nu\in\a^*$ and let $\bar\pi_{\sigma,\nu}$ be the corresponding
principal series representation induced from the parabolic
$\bar{P}=MA\bar N$ opposite to $P$. The representation is defined
to be the right regular representation on the Hilbert space
$H_{\sigma,\nu}$ of all functions $f:G\ra V_\sigma$ satisfying
$f(ma\bar nx)=a^{\nu+\rho_{_P}}\sigma(m)f(x)$ for $ma\bar n\in
MA\bar N$ and $\int_K\norm{f(k)}^2dk<\infty$, modulo
nullfunctions. The space of smooth vectors
$\bar\pi_{\sigma,\nu}^\infty$ coincides with the set of $f$ which
are smooth on $G$. We especially consider $f$ of the form
$$
f(ma\bar nn)\= a^{\nu-\rho_{_P}}\ph(n)
$$
for $n\in N$, where $\ph\in C_c^\infty(N,V_\sigma)$. This
function, defined on the open Bruhat cell $\bar PP$, extends by
zero to a smooth function on $G$. Let
$U_{\sigma,\nu}^\infty\subset H_{\sigma,\nu}^\infty$ denote the
subset of all $f$ of this form.

Let $H_\sigma^\infty$ be the space of all smooth $f:K\ra
V_\sigma$ with $f(mk)=\sigma(m)f(k)$ for $m\in M$, $k\in K$. For
$\nu\in\a^*$ the function
$$
f_\nu(ma\bar nk)\= a^{\nu+\rho_{_P}}\sigma(m)f(k)
$$
defines an element of $H_{\sigma,\nu}^\infty$ and this attachment
sets up an isomorphism of Frech\'et spaces $H_\sigma^\infty\ra
H_{\sigma,\nu}^\infty$ for any $\nu$. Let $U_\sigma^\infty\subset
H_\sigma^\infty$ be the inverse image of $U_{\sigma,\nu}^\infty$
then this space does not depend on $\nu$, which justifies the
notation.

From \cite{wall2}, Theorem 15.4.1 and Theorem 15.6.1 we take

\begin{theorem}
Let $\eta$ be generic and let $\nu\in\a^*$ with $\Re(\nu)<0$. Then
for any $f\in H_\sigma^\infty$ the integral
$$
J_{\sigma,\nu}(f)\=\int_N \eta(n)^{-1} f_\nu(n)dn
$$
converges and extends to a holomorphic map on $\a^*$.

Let $\psi$ be any Whittaker functional on
$\bar\pi_{\sigma,\nu}^\infty$, then
$\psi(U_{\sigma,\nu}^\infty)=0$ implies $\psi=0$. Further, for
any Whittaker functional $\psi$ there is a functional $\mu$ on
$V_\sigma$ such that
$$
\psi(f)\=\mu\left(J_{\sigma,\nu}(f)\right).
$$
\end{theorem}

Fix $\nu$ with $\Re(\nu)<0$. Let $f\in U_\sigma^\infty$ and let
$\ph\in C_c^\infty(N,V_\sigma)$ be the function such that
$f_\nu(ma\bar nn)=a^{\nu+\rho_{_P}}\sigma(m)\ph(n)$. Let $\psi$
be a Whittaker functional and $\mu$ be the corresponding
functional on $V_\sigma$. Then we have
\eqn
W_{f_\nu}(am)&=& \mu\left( \int_N\eta(n)^{-1} f_\nu(nam) dn\right)\\
&=& \int_N\eta(n)^{-1}
a^{\nu+\rho_{_P}}\mu(\sigma(m)\ph(n^{am}))dn.
\endeqn
Thus for $\la\in\a^*$ with $\Re(\la)>k_0\rop$ we get that
$I_{f_\nu}(\xi_\tau,\la)$ equals
$$
\tilde I_f(\xi_\tau,\la +\nu)\=\int_A\int_N\int_M\eta(n)^{-1}
a^{\la +\nu}\mu(\sigma(m)\ph(n^{am}))\overline{\xi_\tau(m)} dm dn
da.
$$
Recall that $\Delta\subset\Phi^+$ denotes the set of simple roots.
Let $N_0=\exp(\bigoplus_{\alpha\notin \N\Delta}\n_\alpha)$ and
$N_\alpha=\exp(\n_\alpha\oplus \n_{2\alpha})$ for $\alpha\in
\Delta$. For $n\in N$ define
$$
\tilde{\ph}(n)\=\tilde\ph_{\mu,\sigma,\tau}(n)\=\int_M\mu(\sigma(m)\ph(n^m))
\overline{\xi_\tau(m)} dm.
$$
We have a canonical identification $C_c^\infty(N,V_\sigma)\cong
C_c^\infty(N)\otimes V_\sigma$. We equip $C_c^\infty(N)$ with the
usual inductive limit topology, then the space $C_c^\infty(N)$
contains as a dense subspace the algebraic tensor product
$$
T\= \left( \bigotimes_{\alpha\in\Delta}C_c^\infty(N_\alpha)\right)
\otimes C_c^\infty(N_0).
$$
Let $T_{\sigma,\nu}\subset U_{\sigma,\nu}$ be the subspace of all
$f$ as above with $\ph=f_\nu|_N$ in $T\otimes V_\sigma$. Now
suppose $\ph$ lies in $T\otimes V_\sigma$, then $\tilde\ph$ lies
in $T$, since the spaces $\n_\alpha$, $\alpha\in\Phi^+$ are
stable under $M$. So suppose that $\tilde\ph$ is a finite sum
$$
\tilde\ph\=\sum_{i=1}^k\left(\prod_{\alpha\in\Delta}\tilde\ph_{i,\alpha}\right)
\tilde\ph_{i,0},
$$
where $\tilde\ph_{i,\alpha}\in C_c^\infty(N_\alpha)$ and
$\tilde\ph_{i,0}\in C_c^\infty(N_0)$. Then
\eqn
I_{f_\nu}(\xi_\tau,\la)&=& \int_A\int_N \eta(n)^{-1}
a^{\la+\nu}\tilde\ph(n^a)dnda\\
&=& \sum_{i=1}^k \int_A a^{\la+\nu}\prod_{\alpha\in\Delta}\int_{N_\alpha}\eta(n_\alpha^{-1})\tilde\ph_{i,\alpha}(n_\alpha^a)
dn_\alpha\ \int_{N_0}\tilde\ph_{i,0}(n^a)dn da
\endeqn
Let $\rop'=\rop-\rez 2\sum_{\alpha\in\N\Delta}m_\alpha\alpha$,
then
$$
\int_{N_0}\tilde\ph_{i,0}(n^a)dn\=
a^{2\rop'}\int_{N_0}\tilde\ph_{i,0}(n)dn.
$$
We may assume that
$$
\int_{N_0}\tilde\ph_{i,0}(n)dn\= 1
$$
for any $i$. Then we get
$$
I_{f_\nu}(\xi_\tau,\la)\=\sum_{i=1}^k \int_A
a^{\la+\nu+2\rop'}\prod_{\alpha\in\Delta}\int_{N_\alpha}\eta(n_\alpha^{-1})\tilde\ph_{i,\alpha}(n_\alpha^a)
dn_\alpha\ da.
$$
We can write $A=\prod_{\alpha\in\Delta}A_\alpha$, where $A_\alpha
=\{ a\in A | a^\beta=1\ {\rm for\ all\ } \beta\in\Delta,\
\beta\ne\alpha\}$. Then we get
$$
I_{f_\nu}(\xi_\tau,\la)\=\sum_{i=1}^k
\prod_{\alpha\in\Delta}\int_{A_\alpha}
a^{\la+\nu+2\rop'}\int_{N_\alpha}\eta(n_\alpha^{-1})\tilde\ph_{i,\alpha}(n_\alpha^a)
dn_\alpha\ da.
$$
Fix an index $i$ and a simple root $\alpha$. Let $G_\alpha$ be
the connected subgroup of $G$ corresponding to the Lie subalgebra
generated by
$$
\n_{-\alpha}\oplus\n_\alpha.
$$
Then $G$ is a real rank one semisimple group with split torus
$A_\alpha$. Let $P_\alpha=P\cap G_\alpha$ and $M_\alpha=M\cap
G_\alpha$, then $P_\alpha=M_\alpha A_\alpha N_\alpha$ is a
minimal parabolic of $G_\alpha$. Mapping a function $h$ on $M$ to
$m_0\mapsto\int_Mh(m_0m)\overline{\xi_\tau(m)}dm$ is an
$L^2$-projection onto the part of $L^2(M)$ spanned by $\xi_\tau$.
Therefore it follows that for any $n\in N_\alpha$:
$$
\int_{M_\alpha}\tilde\ph_{i,\alpha}(n^m)\overline{\xi_{\tau}(m)}
dm\= \tilde\ph_{i,\alpha}(n).
$$
Writing $\xi_\tau$ also for $\xi_{\tau}|_{M_\alpha}$ we get

\begin{lemma}\label{decomp} We have
$$
I_{f_\nu}(\xi_{\tau},\la)\= \sum_{i=1}^k \prod_{\alpha\in\Delta}
I_{f_{i,\alpha,\nu}}^{G_\alpha}(\xi_{\tau},\la+2\rop'),
$$
where $I_{f_{i,\alpha,\nu}}^{G_\alpha}(\xi_{\tau},\la+2\rop')$
denotes the Mellin transform with respect to the group $G_\alpha$
and $f_{i,\alpha,\nu}\in {\rm
Ind}_{P_\alpha}^{G_\alpha}(\sigma\otimes\nu\otimes 1)$ is given by
$$
f_{i,\alpha,\nu}(m_\alpha a_\alpha n_\alpha \theta(n_\alpha)) \=
a_\alpha^{\nu+\rop^\alpha}\sigma(m_\alpha)\ph_{i,\alpha}(n_\alpha),
$$
and
$$
\ph_{i,\alpha}(n_\alpha)\=\int_M\sigma(m)\ph(n_\alpha^m)\overline{\xi_{\tau}(m)}
dm.
$$
\end{lemma}

Note that $I_{f_{i,\alpha,\nu}}^{G_\alpha}(\xi_{\tau},\la+2\rop')$
only depends on the restriction of $\la$ to $A_\alpha$. By the
results of section \ref{rank-one} it follows that
$I_{f_{i,\alpha,\nu}}^{G_\alpha}(\xi_{\tau},\la+2\rop')$ extends
to a meromorphic function and so then does
$I_{f_\nu}(\xi_\tau,\la)$. The position of possible poles can be
read off from Theorem \ref{2.7}

We next want to show that $I_{f_\nu}(\xi_\tau,\la)$ extends to a
meromorphic function for any $f_\nu\in H_{\sigma,\nu}^\infty$. We
need the

\begin{lemma} \label{3.3}
Let $f\in H_{\sigma,\nu}^\infty$, then for every $d\in\N$ there
is a sequence $f_j\in T_{\sigma,\nu}^\infty$ such that
$I_{Df_j}(\xi_\tau,\la)$ converges to $I_{Df}(\xi_\tau,\la)$
locally uniformly on $\Re(\la)>k_0\rop$ for every $D\in U(\g)$ of
degree $\le d$.
\end{lemma}

\prf Fix $d\in\N$. We show the lemma in two steps.  First assume
that $f\in U_{\sigma,\nu}^\infty$, i.e. $f(ma\bar
nn)=a^{\nu-\rop}\sigma(m)\ph(n)$, where $\ph\in
C_c^\infty(N,V_\sigma)$. Then there is a sequence
$(\ph_j)_{j\in\N}$ in $T\otimes V_\sigma$ such that $\ph_j$
converges to $\ph$ in the inductive limit topology. Let $f_j$ be
defined by $f_j(ma\bar nn)=a^{\nu-\rop}\sigma(m)\ph_j(n)$. Let
$N\in\N$, then by Corollary \ref{1.4}  we know we know that there
are $D_1,\dots D_m\in U(\g)$ with
$$
|W_{f_j}(am)|\ \le\ \min_{\alpha\in\Delta} a^{-N\alpha}
\left(\sum_{w\in W}a^{k_0 w\rop}\right)\left(\sum_{k=1}^m
\norm{D_k f_j}\right).
$$
Since $\ph_j$ converges to $\ph$ we conclude that $\norm{D_k
f_j}$ converges to $\norm{D_k f}$ for every $k$, which implies
that there is a constant $C>0$ such that
$$
|W_{f_j}(am)|\ \le\ C\min_{\alpha\in\Delta} a^{-N\alpha}
\left(\sum_{w\in W}a^{k_0 w\rop}\right),
$$
for all $j\in\N$. The claim follows by dominated convergence.

To prove the general case let now $f\in H_{\sigma,\nu}^\infty$
and $d\in\N$. By the first part of this proof it now suffices to
show that there is a sequence $f_j\in U_{\sigma,\nu}^\infty$ such
that $I_{Df_j}(\xi_{\tau},\la)$ converges to
$I_{Df}(\xi_{\tau},\la)$ for every $D\in U(\g)$ of degree $\le d$.
For this let $(\delta_j)_j$ be a sequence of functions $N\ra
[0,1]$ such that
\begin{itemize}
\item
$\delta_j\in C_c^\infty(N)$ for every $j\in\N$, and
\item
$\delta_{j+1}\ge\delta_j$ for every $j\in\N$, and
\item
for every $n\in N$ there is a $j\in\N$ with $\delta_j(n)=1$.
\end{itemize}
Now set $f_j(ma\bar nn)=a^{\nu-\rop}\sigma(m)f(n)\delta_j(n)$,
then $f_j\in U_{\sigma,\nu}^\infty$ and $f_j$ converges pointwise
to $f$.

For $f\in H_{\sigma,\nu}^\infty$ let
$$
S(f)\=\int_N\norm{f(n)}dn.
$$
On p. 382 of \cite{wall2} it is shown that this integral
converges. Then $S$ defines a norm on $H_{\sigma,\nu}^\infty$ and
$$
\norm{J_{\sigma,\nu}(f)}\ \le\ S(f).
$$
We compute
$$
\int_N\norm{f(nam)} dn\= a^{\Re(\nu)+\rop}\int_N\norm{f(n)}dn,
$$
so that
$$
W_f(am)\ \le\ a^{\Re(\nu)+\rho} S(f).
$$
The proof of Proposition \ref{1.2} now gives

\begin{lemma}
There are $D_1,\dots,D_m\in U(\n)$ such that
$$
|W_{f_j}(am)|\ \le\ \min_{\alpha\in\Delta}
a^{-N\alpha}a^{\Re(\nu)+\rop}\sum_{l=1}^m S(D_l f_j).
$$
\end{lemma}

Finally the sequence $\delta_j$ can be chosen so that their
$U(\n)$-derivatives remain bounded, so there is a constant $C>0$
so that $S(D_l f_j)\le C$ for all $l$ and all $j$. This implies
Lemma \ref{3.3}.
\qed

We want to apply the preceding lemma to derive the meromorphic
continuation of $I_{f_\nu}(\xi,\la)$ for general $f_\nu$. In
order to do this we need

\begin{lemma}
Let $f_\nu\in T_{\sigma,\nu}^\infty$ and write as in Lemma
\ref{decomp}:
$$
I_{f_\nu}(\xi_\tau,\la)\= \sum_{i=1}^k\prod_{\alpha\in\Delta}
I_{f_{i,\alpha,\nu}}^{G_\alpha}(\xi_{\tau},\la+2\rop').
$$
For $\alpha\in\Delta$ let $D_\alpha\in U(\Lie (G_\alpha))$, then
there is $D\in U(\g)$ such that
$$
I_{Df_\nu}(\xi_\tau,\la)\= \sum_{i=1}^k\prod_{\alpha\in\Delta}
I_{D_\alpha f_{i,\alpha,\nu}}^{G_\alpha} (\xi_{\tau},\la+2\rop').
$$
\end{lemma}

\prf It suffices to consider the case $D_\alpha=Id$ if
$\alpha\ne\alpha_0$ for some fixed $\alpha_0$ and
$D_{\alpha_0}=X\in \Lie(G_{\alpha_0})$. For this it suffices to
assume that $X$ lies in a generating set of $\Lie(G_{\alpha_0})$,
so, one may take $X$ inside $\n_{-\alpha_0}$ or inside
$\n_{\alpha_0}$. In both cases we take $D=X$ and the claim
follows since the factors associated to $\alpha\ne\alpha_0$ are
annihilated by $X$.
\qed

The meromorphic continuation of $I_{f_\nu}(\xi_\tau,\la)$ now
follows from Lemma \ref{3.3}, Lemma \ref{decomp} and Proposition
\ref{2.6}. This all works for $\Re(\nu)>0$. But since
$$
I_{f_\nu}(\xi_\tau,\la)\=\tilde I_{f}(\xi_\tau+\nu),
$$
we automatically get the meromorphicity for all $\nu$.

We intend to generalize this result to an arbitrary Banach
representation $\pi$. For this we need

\begin{lemma}
Let $\pi$ be an admissible irreducible Banach representation of
$G$, then there are $\sigma\in\hat M$ and $\nu\in\a^*$ such that
the Frech\'et representation $\pi^\infty$ is a quotient of
$\pi_{\sigma,\nu}^\infty$.
\end{lemma}

\prf Let $\pi_K$ denote the admissible irreducible $(\g,K)$-module
of $K$-finite vectors in $\pi$. The dual $(\g,K)$-module $\breve
\pi_K$ also is admissible and irreducible. By the
subrepresentation theorem there exist $\sigma\in\hat M$ and
$\nu\in\a^*$ such that $\breve \pi_K$ injects into
$\bar\pi_{\breve \sigma,\nu,K}$. Dualizing we get a nontrivial
and hence surjective $(\g,K)$-map from $\bar\pi_{\sigma,\nu,K}$ to
$\pi_K$, which proves the claim for the underlying
$(\g,K)$-modules. By Corollary 10.5 of\cite{cas} the claim
follows.
\qed

So let finally $\pi$ be an arbitrary admissible irreducible
Banach representation of $G$. Let $\phi\in\pi^\infty$ and let
$F:\bar\pi_{\sigma,\nu}^\infty\ra\pi^\infty$ be a surjective
homomorphism. Pick some $f\in\bar\pi_{\sigma,\nu}^\infty$ such
that $F(f)=\phi$. Then $\psi\circ F$ is a Whittaker functional on
$\bar\pi_{\sigma,\nu}^\infty$. It then follows
$$
I_\phi(\xi_\tau,\la)\= I_f(\xi_\tau,\la),
$$
so the meromorphic continuation of $I_\phi(\xi_\tau,\la)$ is
established. Theorem \ref{higher-rank} is proven.

University of Exeter\\
Mathematics\\
Exeter EX4 4QE\\
Devon, UK

%\klghlkdshg

\newpage
\begin{thebibliography}{XXX}

\bibitem{CLP-S}
 \bf Cogdell, J.; Li, J.; Pyatetskii-Shapiro, I.:
 \it The meromorphic continuation of Kloosterman-Selberg Zeta functions.
 \rm Springer Lecture Notes 1422, 23-35 (1990).

\bibitem{CLP-SS}
 \bf Cogdell, J.; Li, J.; Pyatetskii-Shapiro, I.; Sarnak, P.:
 \it Poincar\'e series for $SO(n,1)$.
 \rm Acta Math. 167, 229-285 (1991).

\bibitem{bump}
 \bf Bump, D.:
 \it Automorphic Forms on $GL(3,\R)$.
 \rm Lecture Notes in Mathematics 1083. Springer-Verlag 1984.


\bibitem{cas}
 \bf Casselman, W.:
 \it Canonical extensions of Harish-Chandra modules to representations of G.
 \rm Can. J. Math. 41, 385-438 (1989).

\bibitem{DM}
 \bf Dixmier, J.; Malliavin, P.:
 \it Factorisations de fonctions et de vecteurs indefiniment differentiables.
 \rm Bull. Sci. Math., II. Ser. 102, 305-330 (1978).

\bibitem{goldf}
 \bf Goldfeld, D.:
 \it Kloosterman Zeta Functions for $GL(n,\Z)$.
 \rm Proceedings of the ICM, Berkeley, CA 1986.

\bibitem{humpliealg}
 \bf Humphreys, J.:
 \it Introduction to Lie Algebras and Representation Theory.
 \rm Springer-Verlag 1972.

\bibitem{rag}
 \bf Raghunathan, M.:
 \it Lattices in Lie Groups
 \rm New York. Springer-Verlag. 1972.

\bibitem{Wall-asymp}
 \bf Wallach, N.:
 \it Asymptotic expansions of generalized matrix coefficients of real reductive groups.
 \rm Lie Group Representations I. SLNM 1024, 287-369 (1983)

\bibitem{wall2}
 \bf Wallach, N.:
 \it Real Reductive Groups II.
 \rm Academic Press 1993.

\end{thebibliography}
\end{document}